\documentclass[review,hidelinks,onefignum,onetabnum]{siamart220329}

\usepackage{mathrsfs}
\usepackage{eepic}
\usepackage{amsfonts}
\usepackage{verbatim}
\usepackage{ragged2e}
\usepackage[numbers,sort&compress]{natbib}
\usepackage{arydshln}

\usepackage{url}
\usepackage{color}
\usepackage{booktabs}
\usepackage{threeparttable}
\usepackage{amssymb}
\usepackage{bbding}
\usepackage{booktabs}
\usepackage{graphicx}
\usepackage{dsfont}
\usepackage{amsmath}
\usepackage{float}
\usepackage{subfigure}
\usepackage{multirow}
\usepackage{multicol}
\usepackage[figureright]{rotating}
\usepackage{lscape}
\usepackage{epstopdf}

\usepackage{lipsum}
\usepackage{amsfonts}
\usepackage{graphicx}
\usepackage{epstopdf}
\usepackage{algorithmic}
\ifpdf
  \DeclareGraphicsExtensions{.eps,.pdf,.png,.jpg}
\else
  \DeclareGraphicsExtensions{.eps}
\fi

\newsiamremark{remark}{Remark}
\newsiamremark{hypothesis}{Hypothesis}
\crefname{hypothesis}{Hypothesis}{Hypotheses}
\newsiamthm{claim}{Claim}

\headers{A DR index tracking model with the CVaR penalty}{Ruyu Wang, Yaozhong Hu, Chao Zhang}

\title{A distributionally robust index tracking model with the CVaR penalty: tractable reformulation\thanks{Submitted to the editors DATE.
\funding{The work is supported by ``the Natural Science Foundation of China" (Grant No. 12171027). Ruyu Wang is supported by China Scholarship Council.
}}}

\author{Ruyu Wang\thanks{School of Mathematics and Statistics, Beijing Jiaotong University
(\email{wangruyu@bjtu.edu.cn}, Correspondence:\email{zc.njtu@163.com}). }
\and Yaozhong Hu\thanks{Department of Mathematical and Statistical Sciences, University of Alberta	(\email{yaozhong@ualberta. ca}). }
\and Chao Zhang\footnotemark[2]}

\usepackage{amsopn}

\ifpdf
\hypersetup{
	pdftitle={A distributionally robust index tracking model with the CVaR penalty: tractable reformulation},
	pdfauthor={Ruyu Wang, Yaozhong Hu, Chao Zhang}
}
\fi

\begin{document}
	
	\maketitle
	
	\begin{abstract}
We propose a distributionally robust index tracking model with the conditional value-at-risk (CVaR) penalty. The model combines the idea of distributionally robust optimization for data uncertainty and the CVaR penalty to avoid large tracking errors. The probability ambiguity is described through a confidence region based on the first-order and second-order moments of the random vector involved. We reformulate the model in the form of a min-max-min optimization into an equivalent nonsmooth minimization problem. We further give an approximate discretization scheme of the possible continuous random vector of the nonsmooth minimization problem, whose objective function involves the maximum of numerous but finite nonsmooth functions. The convergence of the discretization scheme to the equivalent nonsmooth reformulation is shown under mild conditions. A smoothing projected gradient (SPG) method is employed to solve the discretization scheme. Any accumulation point is shown to be a global minimizer of the discretization scheme. Numerical results on the NASDAQ index dataset from January 2008 to July 2023 demonstrate the effectiveness of our proposed model and the efficiency of the SPG method, compared with several state-of-the-art models and corresponding methods for solving them.
	\end{abstract}
	
	\begin{keywords}
	{Index tracking}, Distributionally robust optimization, Conditional value-at-risk, 
	{Nonsmooth}, {Smoothing projected gradient method}.
	\end{keywords}
	
	\begin{MSCcodes}
		46N10, 
		57R12, 
		90C15, 
		90C47
	\end{MSCcodes}
	
\section{Introduction}\label{sec1}
		
		Index tracking aims to replicate the index of a financial market by constructing a portfolio consisting of assets in that market that minimizes the tracking error, which measures how closely the portfolio mimics the performance of the benchmark. Index tracking has received a great deal of attention in the literature. Xu et al. \cite{xu2016efficient} consider an index tracking model that minimizes a quadratic tracking error while enforcing an upper bound on the number of assets in the portfolio. They also propose a non-monotone projected gradient method for solving the model. Zhang et al. \cite{zhang2019robust} propose a robust and sparse index tracking model, by adding the $\ell_{2}$ and $\ell_{p}$ ($0 < p < 1$) norm penalties in the objective function. The SPG method is employed to solve it \cite{zhang2019robust}. {We emphasize that the original index model ignores extreme risk control. This means a portfolio may suffer from market jump risk because it closely replicates a benchmark index. Hence, it is crucial to control downside risk. Jorion \cite{jorion2003portfolio} adds a ``Value at Risk" (VaR) penalty in the objective function to the index tracking model, which is a risk measure that controls downside risk.} But the VaR has some undesirable characteristics such as a lack of subadditivity and convexity. As an alternative measure of risk, the CVaR is known to have better properties than the VaR; see Rockafellar et al. \cite{rockafellar2000optimization} for details of the CVaR. Goel et al. \cite{goel2018index} propose to use the generalized forms of the CVaR for an index tracking model. Wang et al. \cite{wang2012mixed} introduce a CVaR risk constraint into the index tracking model to control the downside risk of tracking portfolios. The CVaR constraint helps to prevent investors from large losses. {When a CVaR constraint is added to the general index tracking model, the resulting index tracking problem can be transformed into a mixed 0-1 linear programming (LP) problem. When the number of 0–1 variables is relatively small, the authors show that the resulting mixed 0–1 linear program can be efficiently solved using standard optimization software. However, implementing the CVaR constraint needs defining the investor's maximum risk tolerance, a value that is challenging to accurately estimate.} In this paper, we consider minimizing the objective function consisting of both the tracking error and the risk measured using the CVaR.

		Recently, there has been extensive use of distributionally robust optimization (DRO) in real applications including index tracking and data mining \cite{kang2022robust,delage2010distributionally,wiesemann2014distributionally,guo2019distributionally,li2020fast,bertsimas2019adaptive}. For example, Kang et al. \cite{kang2022robust} formulate a DRO index tracking model by integrating uncertain information on the distributions specified by the $\phi$-divergence. Liu et al. \cite{liu2019closed} come up with the closed-form optimal solutions for the distributionally robust mean-CVaR portfolio selection problems. Because the objective function of this problem is a simple linear function, solving the maximization problem in robust optimization yields a closed-form solution. However, due to the non-linear nature of the objective function in the more general index tracking problems, it is difficult to get the closed-form solution of the index tracking problem. Xu et al. in \cite{xu2018distributionally} analyze the condition for the strong Lagrange duality of the minimax DRO problems in detail. They also propose a discretization scheme for solving the dualized DRO. Chen et al. \cite{chen2019discrete} broaden the discretization scheme to solve a two-stage distributionally robust linear complementarity problem. Jiang and Chen \cite{jiang2023pure} propose an approximation problem with discretization to solve distributionally robust pure characteristics demand models with stochastic complementarity constraints. 
		DRO is especially attractive for problems in which input data and their corresponding distributions are uncertain. In practice, the true probability distribution of the random vector can hardly be acquired in index tracking. This motivates us to further enroll the idea of DRO in the index tracking model. 
		
		In this paper, we propose a distributionally robust index tracking model with $\ell_2$ norm penalty and the CVaR penalty in the objective function, which we call DRCVaR index tracking model for short. {The probability ambiguity is described through a confidence region based on the first-order and second-order moments of the random vector involved.} The model is promising, since we combine the idea of distributionally robust to index tracking error that reflects information uncertainty well, the $\ell_2$ norm for avoiding over-fitting, as well as the CVaR penalty for risk control in the original index tracking model. However, these terms also make it highly challenging with tractability and efficiency of computation, because it is a min-max-min optimization problem. By the minimax theorem in \cite{sion1958general}, we equivalent transform DRCVaR into a min-max optimization problem. {Since the Slater condition of the inner maximization is satisfied in the above probability ambiguity set, we equivalently transform the maximization problem by the Lagrange dual to a semi-definite programming (SDP).} Following this result, we establish an equivalent tractable reformulation for DRCVaR using Lagrangian dual. Meanwhile, we prove that the dual problem is convex and the set of the optimal solutions is nonempty. {The objective of the dual problem includes a maximization of an infinite number of convex functions if the underlying distribution is continuous.} We address this problem with a discretization scheme and provide the convergence. So far, we have transformed DRCVaR to the problem involving the maximization of numerous nonsmooth functions.

		
		It is known that the smoothing methods are efficient for dealing with nonsmooth optimization \cite{nesterov2005smooth,zhang2009smoothing,zhang2020smoothing,bian2020smoothing,bian2015linearly,chen2012smoothing,xu2014smoothing,xu2015smoothing}. Nesterov \cite{nesterov2005smooth} approximates the nonsmooth objective function by a function with Lipschitz continuous gradient and minimizes the smooth function by an efficient gradient method. Zhang et al. \cite{zhang2009smoothing} propose the SPG method for solving the minimization problem on a closed convex set. Zhang and Chen \cite{zhang2020smoothing} propose a smoothing active set method for minimizing the linearly constrained nonLipschitz nonconvex problems. Chen et al. \cite{chen2012stochastic} propose a smoothing sample average approximation method to minimize an expected residual function. Bian and Chen \cite{bian2020smoothing} develop a smoothing proximal gradient algorithm for a class of constrained optimization problems with the objective function defined by the sum of a nonsmooth convex function and a cardinality function. The nonsmooth term arising from the dual and each nonsmooth convex function induced from the CVaR penalty can be approximated using a smoothing function. Compared with the nonsmooth term from the dual and each nonsmooth convex function from the CVaR penalty, the maximum of numerous convex functions is more difficult to deal with. Motivated by the above properties of existing index tracking models and the challenges of numerical tractability, we aim to transform the DRCVaR model equivalently and use the SPG method to solve the proposed DRCVaR model efficiently.

		In Section \ref{sec3}, we propose the DRCVaR model, and transform the inner maximization in the ambiguity set to an SDP. In Section \ref{equivalent}, we establish an equivalent tractable reformulation for DRCVaR using Lagrangian dual. Meanwhile, we prove that the dual problem is convexity and the set of the optimal solutions is nonempty. {The objective of the dual problem includes a maximization of an infinite number of convex functions.} We address this problem with a discretization scheme and study the convergence in Section \ref{disscheme}. In Section \ref{sec4}, we use the SPG method to solve problems involving the maximization of numerous nonsmooth functions. The convergence result is also provided. In Section \ref{sec5}, we illustrate the numerical implementation of our method for solving a distributionally robust portfolio model for index tracking. Numerical results demonstrate the merits of the DRCVaR model and the SPG method. We make a conclusion in Section \ref{sec6}.

		Throughout the paper, we use the following notation. We use $\mathbb{S}_+^{d}$ to denote the cone of symmetric positive semi-definite matrices. We denote by $\|x\|$ the Euclidean norm of a vector $x$. For a number $r\in \mathbb{R}$, we denote $[r]_{+} = \max\left\{r, 0\right\}$. Given two square matrices $M$ and $N$, we write $M \preceq N$ to indicate $N - M$ being positive semi-definite. {The notation $\left\langle M,N \right\rangle$ refers to the inner product between two matrices $M$ and $N$ of the same dimension, i.e., $\left\langle M,N \right\rangle= \sum\limits_{i,j}M_{i,j}N_{i,j}$.} We denote by $B(x,\delta):=\left\lbrace y \in \mathbb{R}^{d} :\|x-y\|\leq \delta\right\rbrace $ the closed ball of radius $r$ centered at $x\in \mathbb{R}^d$.

			\section{DRCVaR index tracking model}\label{sec3}
			
			In this section, we establish the DRCVaR model for index tracking, provide the motivation of this model, and analyze the properties of the model. {We briefly review some important index tracking models as follows.
				\subsection{Deterministic models using exactly historical data}\label{secmodel}}
			
			
			Various index portfolio models have been proposed to give an optimal tracking portfolio for the investor, based directly on historical data. Then, the optimal tracking portfolio obtained from in-sample data (historical data) will be used to track the index in the out-of-samples.

			Let $\widehat{\xi_{B}^j} =(\widehat{\xi_{B,1}^j},\ldots,\widehat{\xi_{B,d}^j})\in \mathbb{R}^{d}$ be the 
			observed historical individual 
			return vector of the $d$ assets, and {$\widehat{\xi_{a}^j}\in \mathbb{R}$} be the 
			observed corresponding 
			random market index return for $j=1,\ldots,N$. We denote by $\widehat{\xi^j} = (\widehat{\xi_B}^T,\widehat{\xi_a^j})^T \in \mathbb{R}^{d+1}$. Let $x = (x_1,\ldots,x_d)^T \in \mathbb{R}^d$ be the tracking portfolio, with $x_i$ being the investment weight in the $i$th component stock. 
			%
			%
			Various portfolio index tracking models can be outlined in the following model, using historical data. 
			\begin{eqnarray}\label{general}
				\min\limits_{x\in X}~\frac{1}{N} \sum\limits_{j=1}^N\left[ \psi\left(\widehat{\xi_{a}^j} -\widehat{\xi_{B}^j}^{\top} x\right)\right] + \tau R(x).
			\end{eqnarray}
			Here the first term is the data fidelity term, the convex function $\psi$ tries to describe the deviation between the return rate of the tracking portfolio and the return rate of the index. {For instance, we can choose $\psi(a)=\vert a\vert$, or $\psi(a)=a^2$}. The function $R$ in the second term is the penalization term that can help to control the risk, reflect the sparsity of the portfolio, and enhance the out-of-sample performance, etc. The parameter $\tau$ balances the trade-off of the data fidelity term and the penalty term. The feasible set $X$ is a convex set, which includes some prior knowledge of the desired portfolio, e.g., the sum-to-one requirement $\sum\limits_{i=1}^d x_i =1$. 
			
			{We outline some index tracking models that can be written in the general model \eqref{general} as follows.}
			
			-- The Lasso sparse index tracking problem in \cite{sant2020lasso} can be formulated as
			\begin{eqnarray*}
				(\mathrm{Lasso~sparse})\qquad
				\begin{aligned}
					\min\limits_{x} 	& ~\frac{1}{N} \sum\limits_{j=1}^N\left(\widehat{\xi_{a}^j} -\widehat{\xi_{B}^j}^{\top} x\right)^2+\tau_1\|x\|_1 \\
					\text { s.t. }& \sum\limits_{i=1}^d x_i =1,
				\end{aligned}
			\end{eqnarray*}
			{where $\tau_1$ is a given regularization parameter, and the $\ell_1$-norm penalty aims to enhance the sparsity of the portfolio.} It's worth noting that there is no nonnegativity constraint. Therefore, negative weights may appear, which will lead to the possibility of short selling. The $(\mathrm{Lasso~sparse})$ model is solved using the optimization software ``CPLEX".

			{After introducing the downside risk control with the CVaR risk measure in \cite{wang2012mixed}, the mixed 0-1 LP index tracking model with the CVaR risk constraint is given as follows:}
			\begin{eqnarray*}
				(\mathrm{Mixed~0-1~LP})\qquad
				\begin{array}{ll}
					\min\limits_{x, Z,\alpha} & \frac{1}{N} \sum\limits_{j=1}^N \left|\widehat{\xi_{a}^j} -\widehat{\xi_{B}^j}^{\top} x\right| \\
					\text { s.t. } & \sum\limits_{i=1}^d Z_i=K \\
					& l_i Z_i \leq x_i \leq u_i Z_i, \quad i=1,2, \ldots, d \\
					& \sum\limits_{i=1}^d x_i=1 \\
					& Z_i \in\{0,1\}, \quad i=1,2, \ldots, d\\
					&\alpha+\frac{1}{(1-\beta)N} \sum\limits_{j=1}^N p_j\leq \eta,\\
					&p_j\geq -\widehat{\xi_{B}^j}^{\top} x-\alpha,\quad p_j\geq 0, \quad j=1,2, \ldots, N,
				\end{array}
			\end{eqnarray*}
			where $K$ is the number of stocks selected to track the index, $Z_i=1$ means that the $i$th stock is included in the tracking portfolio, $Z_i=0$ means the opposite for $i=1,\ldots,d$. The constants $l_i$ and $u_i$ are the lower and upper bounds of the investment proportion on asset $i$, and $0 < l_i < u_i < 1$. Constraints $l_i Z_i \leq x_i \leq u_i Z_i$, $i=1,2, \ldots, d$ show that if asset $i$ is not selected in the tracking portfolio (i.e., $Z_i = 0$), then $x_i = 0$, and if asset $i$ is selected in the tracking portfolio (i.e., $Z_i = 1$), hence the value of $x_i$ is limited in the interval $[l_i,u_i]$. The $(\mathrm{Mixed~0-1~LP})$ model is solved using the optimization software ``CPLEX".

			-- 
			{Xu et al. \cite{xu2016efficient} consider a model for index tracking, which minimizes a quadratic tracking error while enforcing an upper bound on the number of assets in the portfolio:}
			\begin{eqnarray*}
				(\mathrm{TE-\ell_0})\qquad
				\begin{array}{ll}
					\min\limits_{x} & ~ \frac{1}{N} \sum\limits_{j=1}^N \left(\widehat{\xi_{a}^j} -\widehat{\xi_{B}^j}^{\top} x\right)^2\\
					\text { s.t. } &  \sum\limits_{i=1}^d x_i=1, \\
					& \Vert x\Vert_0\leq K,\\
					&0\leq x_i\leq u, \quad i=1,2, \ldots, d,
				\end{array}
			\end{eqnarray*}
			where $\Vert x\Vert_0$ denotes the number of nonzero entries of $x$, $K$ is a given positive integer, and $u\in [1/K, 1]$ is an upper bound on the weight of each index constituent. In \cite{xu2016efficient}, Xu et al. propose a nonmonotone projected gradient method for solving this problem.


			%
			%
			%

			-- Zhang et al. \cite{zhang2019robust} consider the following robust and sparse nonconvex index tracking model:
			\begin{eqnarray*}
				(\ell_2-\ell_{p})\qquad
				\begin{array}{ll}
					\min\limits_{x} & ~ \frac{1}{N} \sum\limits_{j=1}^N \left(\widehat{\xi_{a}^j} -\widehat{\xi_{B}^j}^{\top} x\right)^2 + \tau_{1}\|x\|^{2} + \tau_{2}\|x\|^{p}_{p},\\
					\text { s.t. } & \sum\limits_{i=1}^d x_i=1,\quad x\geq 0,
				\end{array}
			\end{eqnarray*}
			where $0<p<1$, $\tau_{1},~\tau_{2} \geq 0$ are given parameters. The $\ell_{p}$ norm $(0<p<1)$ has a theoretical guarantee of sparseness, {which is powerful to obtain a sparse portfolio and consequently helps to decrease the transaction fee. The $\ell_{2}$-norm penalty aims to increase the out-of-sample performance. The shortsale (nonnegativity) constraint $x\geq 0$ requires that selling assets short is forbidden, i.e., the fund managers cannot sell the assets that they do not own currently. The shortsale constraint in portfolio models may help to stabilize the stock market \cite{engelberg2018short,chen2019short,anufriev2013impact,grullon2015real}.} Zhang et al. in \cite{zhang2019robust} propose a SPG algorithm to solve it.

			{\subsection{DRCVaR model}}
			
			It is known that the performance of a portfolio on the out-of-samples is more important than that on the in-samples. Since the portfolio will be used for prediction later when the true return rate of the index and the true return rate of the portfolio is not known yet, it is often the case that the in-sample performance is very good by solving a certain index tracking model, however, the out-of-sample performance is not as good as expected. To overcome this difficulty, instead of using historical data exactly in index tracking model, {we consider ${\xi_{B}} =({\xi_{B,1}},\ldots,{\xi_{B,d}})\in \mathbb{R}^{d}$ to be the random return vector and ${\xi_{a}}\in \mathbb{R}$ to be the corresponding random market index return. Let the random vector $\xi: = (\xi_B^T, \xi_a)^T \in \mathbb{R}^{d+1}$. We assume that $\xi$ is distributed according to the probability distribution $P$. One may extend \eqref{general} to the following stochastic index tracking portfolio model with the CVaR penalty:}
			\begin{eqnarray*}
				(\mathrm{SCVaR-P})\qquad
				\min\limits_{x\in X}~ \mathbb{E}_{P}[\psi(\xi_{a} -\xi_{B}^{\top}x)] +\tau_{1}\|x\|^{2}+\tau_{2} \phi_\beta(x),
			\end{eqnarray*}
			where $x = (x_{1}, \ldots, x_{d})^{\top}$ is the decision vector, with each component indicating the proportion of the total amount of money invested in asset $i \in d$,  
			and the regularization parameters $\tau_{1},~\tau_{2}\geq 0$. Here, the CVaR term
			\begin{eqnarray}\label{cvar}
				\phi_\beta(x)= \min\limits_{\alpha\in\mathbb{R}} \left\{\alpha+(1-\beta)^{-1}\mathbb{E}_{P}\left[\ell(x,\xi_{B})-\alpha\right]_{+}\right\}.
			\end{eqnarray}
			Here $\ell(x,\xi_{B})$ is a convex loss function in $x\in X$. In this paper, the convex loss function
			\begin{eqnarray*}
				\ell(x,\xi_{B}):=-x^{\top}\xi_{B}
			\end{eqnarray*}
			is the same as \cite{rockafellar2000optimization}. The $(\mathrm{SCVaR-P})$ model can be solved using the stochastic subgradient (S-Subgrad) method as described in \cite{nemirovski2009robust,wang2022stochastic}. 

			\vskip 2mm
			
			{In practice, it is reasonable to assume that $\xi$ follows the true probability distribution $P$, which is not known exactly but is sure to be contained in an ambiguity set $\mathscr{P}$ constructed using partial information such as historical data.} This motivates us to further enroll the idea of DRO in the index tracking model. We therefore propose in this paper the DRCVaR model formulated as follows:
			\begin{eqnarray}\label{orip}
				(\mathrm{DRCVaR})\qquad
				\min\limits_{x\in \Delta_{d}}\max\limits_{P\in\mathscr{P}}~ \mathbb{E}_{P}[\psi(\xi_{a} -\xi_{B}^{\top}x)]
				+\tau_1 \|x\|^2 + \tau_2 
				\phi_{\beta}(x),
			\end{eqnarray}
			where $\Delta_{d}=\left\{x\mid\sum\limits_{i=1}^d x_i=1,~ x\geq 0\right\}$. Let us denote by $\hat{\mu}$ and $\hat{\Sigma}$ the reference values of the mean vector and covariance matrix of the historical data $(\widehat{\xi}^1,\ldots,\widehat{\xi}^N)$. We assume that $\hat{\Sigma}$ is a symmetric and positive definite matrix. The ambiguity set $\mathscr{P}$ in (\ref{orip}) is constructed by moments constraints as intensively studied by Delage and Ye in \cite{delage2010distributionally} and by Xu et al. \cite{xu2018distributionally}:
			\begin{eqnarray}\label{ambiguityset}
				\mathscr{P}=
				\left\{ \begin{array}{l|l}
					P\in \mathcal{M} & \begin{array}{l}
						\left(\mathbb{E}_{P}[\xi]-\hat{\mu}\right)^{\top}\hat{\Sigma}^{-1}\left( \mathbb{E}_{P}[\xi]-\hat{\mu}\right) \leq \kappa_{1}\\
						\mathbb{E}_{P}\left[ (\xi-\hat{\mu})(\xi-\hat{\mu})^{\top}\right] \preceq \kappa_{2} \hat{\Sigma}
					\end{array}
				\end{array} \right\},
			\end{eqnarray}
			where $\kappa_{1}>0$, $\kappa_{2}>0$, $\mathcal{M}$ is the convex set of all probability measures in the measurable space $(\Xi, \mathds{B})$, with $\Xi \subseteq \mathbb{R}^{d+1}$ being a convex compact set known to contain the support of $P$, and $\mathds{B}$ being the Borel $\sigma$-algebra on $\Xi$. 
			
		
		\begin{lemma}[Section A.5.5 of \cite{boyd2004convex}]\label{schur}
			Consider a symmetric matrix $\hat{M}$ partitioned as
			$$
			\hat{M}=\left[\begin{array}{cc}
				A & B \\
				B^\top & C
			\end{array}\right],
			$$
			where $A$ is a symmetric matrix. If $det(A) \neq 0$, the matrix
			$$
			\hat{S}=C-B^{\top}A^{-1}B
			$$
			is called the Schur complement of $A$ in $\hat{M}$. If $A\succ 0$, then $\hat{M}\succeq 0$ if and only if $\hat{S}\succeq 0$.
			
		\end{lemma}

		Below, we show the strong duality holds between the inner maximization problem \eqref{max1} and its Lagrangian dual.
		It is easy to observe that the first constraint in \eqref{ambiguityset} by Lemma \ref{schur} can be equivalently written as
		\begin{eqnarray*}
			\mathbb{E}_{P}\left[\begin{aligned}
				&-\hat{\Sigma}&\hat{\mu}-\xi\\
				&(\hat{\mu}-\xi)^\top &-\kappa_{1} 
			\end{aligned}\right] \preceq 0.
		\end{eqnarray*}

		Let $x\in\Delta_{d}$ and $\alpha\in\mathbb{R}$ be fixed. We consider a maximization problem in the ambiguity set \eqref{ambiguityset}
		\begin{eqnarray}\label{max1}
			\begin{aligned}
				\max\limits_{P\in\mathcal{M}}\quad&\mathbb{E}_{P}[\hat{K}(x,\alpha,\xi)]\\
				\rm{s.t.}\quad& 		\mathbb{E}_{P}\left[\begin{aligned}
					&-\hat{\Sigma}&\hat{\mu}-\xi\\
					&(\hat{\mu}-\xi)^\top &-\kappa_{1} 
				\end{aligned} \right] \preceq 0,\\
				\quad&		\mathbb{E}_{P}\left[ (\xi-\hat{\mu})(\xi-\hat{\mu})^{\top}\right] \preceq \kappa_{2} \hat{\Sigma},\\
				\quad&		\mathbb{E}_{P}\left[ 1\right] =1,
			\end{aligned}
		\end{eqnarray}
		where
		{\begin{eqnarray}\label{K}
				\hat{K}(x,\alpha,\xi)=\psi(\xi_{a} -\xi_{B}^{\top}x)+
				\frac{\tau_{2}}{1-\beta}\left[-x^{\top}\xi_{B}-\alpha\right]_{+}+\tau_{1}\|x\|^{2}
				+\tau_{2}\alpha.
			\end{eqnarray} is a convex function.}
		
		{In Proposition \ref{propinn}, we transform the inner maximization problem \eqref{max1} by the Lagrange dual to a semi-infinite program. For DRO, the equivalence between the primal problem and the dual problem has been well established under the circumstance where the moment problem satisfies the Slater type condition. We also show that the Slater condition of \eqref{max1} is satisfied for the probability ambiguity set \eqref{ambiguityset} in Proposition \ref{propinn}.}

		\begin{proposition}\label{propinn}
			For fixed $x\in\Delta_{d}$ and $\alpha\in\mathbb{R}$, suppose that $\kappa_{1}>0$, $\kappa_{2}>0$, and $\hat{\Sigma}\succ 0$. 
			Denote
			\begin{eqnarray}\label{h1}
				\begin{aligned}
					&h_{1}(x,\alpha,q,\Lambda):=\left\langle \kappa_{2}\hat{\Sigma},\Lambda\right\rangle+(\Lambda^{\top}\hat{\mu}+q)^{\top}\hat{\mu}+ \tau_{1}\|x\|^{2}+\tau_{2}\alpha,\\
					&h_{2,\xi}(x,q,\Lambda):=\psi(\xi_{a} -\xi_{B}^{\top}x)-(\Lambda^{\top}\xi+q)^{\top}\xi+\frac{\tau_{2}}{1-\beta}\left[-x^{\top}\xi_{B}-\alpha\right]_{+}.
				\end{aligned}
			\end{eqnarray}
			Then, the equivalence holds between \eqref{max1} and
			\begin{eqnarray}\label{innerdual}
				\begin{aligned}
					\min\limits_{r,q,\Lambda}\quad&r +h_{1}(x,\alpha,q,\Lambda)+\sqrt{\kappa_{1}}\Vert\hat{\Sigma}^{\frac{1}{2}}(q+2\Lambda\hat{\mu})\Vert\\
					\rm{s.t.}\quad& h_{2,\xi}(x,q,\Lambda)\leq r ,~\forall\xi\in\Xi,\\
					\quad&r\in \mathbb{R},\quad q\in\mathbb{R}^{d+1}, \quad\Lambda\in\mathbb{S}_+^{d+1},
				\end{aligned}
			\end{eqnarray}
			where $r\in \mathbb{R}$, $q\in\mathbb{R}^{d+1}$, and $\Lambda\in\mathbb{R}^{(d+1)\times (d+1)}$ are the dual variables.
			Moreover, the set of optimal solutions of the dual problem \eqref{innerdual} is nonempty and bounded.
		\end{proposition}
		\begin{proof}
			Since the definition \eqref{K} of $\hat{K}(x,\alpha,\xi)$, the Lagrange function of \eqref{max1} is
			\begin{eqnarray*}
				&&\int
				\psi(\xi_{a} -\xi_{B}^{\top}x)-(\Lambda^{\top}\xi-2\hat{q})^{\top}\xi+2\hat{\mu}^{\top}\Lambda\xi+\frac{\tau_{2}}{1-\beta}\left[-x^{\top}\xi_{B}-\alpha\right]_{+}-r
				~dP(\xi)\\
				&&\quad+\left\langle \kappa_{2}\hat{\Sigma},\Lambda\right\rangle -(\Lambda^{\top}\hat{\mu}+2\hat{q})^{\top}\hat{\mu}+ \tau_{1}\|x\|^{2}+\tau_{2}\alpha+r+\left\langle \hat{\Sigma},Q \right\rangle+\kappa_{1}s.
			\end{eqnarray*}
			Here $r \in \mathbb{R}$ and $\Lambda \in \mathbb{R}^{(d+1)\times (d+1)}$ are the dual variables for the last and the second constraint of \eqref{max1}, respectively. And $Q \in \mathbb{R}^{(d+1)\times (d+1)}, \hat{q} \in \mathbb{R}^{d+1}$ and $s \in \mathbb{R}$ form together a matrix consisting of the dual variables for the first constraint of \eqref{max1}. Then the Lagrange dual of \eqref{max1} can take the following form
			\begin{eqnarray}
				\min\limits_{r, \Lambda, Q, \hat{q}, s}&& \left\langle \kappa_{2}\hat{\Sigma},\Lambda\right\rangle -(\Lambda^{\top}\hat{\mu}+2\hat{q})^{\top}\hat{\mu}+ \tau_{1}\|x\|^{2}+\tau_{2}\alpha+r+\left\langle \hat{\Sigma},Q \right\rangle+\kappa_{1}s \label{qqobj}\\
				\text { s.t.}\quad&&\psi(\xi_{a} -\xi_{B}^{\top}x)-(\Lambda^{\top}\xi-2\hat{q})^{\top}\xi+2\hat{\mu}^{\top}\Lambda\xi+\frac{\tau_{2}}{1-\beta}\left[-x^{\top}\xi_{B}-\alpha\right]_{+}\nonumber\\
				&&\qquad -r\leq 0,~\forall \xi \in \mathcal{S}, \nonumber\\
				&& r\in \mathbb{R},\quad \Lambda \succeq 0, \nonumber\\
				&& {\left[\begin{array}{cc}
						Q & \hat{q} \\
						\hat{q}^{\top} & s
					\end{array}\right] \succeq 0.\label{qq}}
			\end{eqnarray}
			
			Similar to the procedure in the proof of Lemma 1 in \cite{delage2010distributionally}, we can simplify the dual problem by analytically solving for the variables $(Q, \hat{q}, s)$ while keeping $(\Lambda, r)$ fixed. For the sake of completeness, we briefly outline these steps below. In view of the semi-definite constraint \eqref{qq}, we consider two cases for the variable $s^*$: either $s^*=0$ or $s^*>0$. Let's first consider the case of $s^* = 0$. In this scenario, if $\hat{q}^* \neq 0$, it would lead to $\hat{q}^{* \top} \hat{q}^* > 0$ and
			\begin{eqnarray*}
				\left[\begin{array}{c}
					\hat{q}^* \\
					y
				\end{array}\right]^{\top}\left[\begin{array}{cc}
					Q^* & \hat{q}^* \\
					\hat{q}^{* \top} & s^*
				\end{array}\right]\left[\begin{array}{c}
					\hat{q}^* \\
					y
				\end{array}\right]=\hat{q}^{* \top} Q^* \hat{q}^*-2 \hat{q}^{* \top} \hat{q}^* y<0, ~\text { for }~ y>\frac{\hat{q}^{* \top} Q^* \hat{q}^*}{2 \hat{q}^{* \top} \hat{q}^*},
			\end{eqnarray*}
			which contradicts equation \eqref{qq}. Therefore, it must be the case that $\hat{q}^* = 0$. Furthermore, given that $\hat{\Sigma}\succ0$, $Q\succeq 0$, and the objective is to minimize $\left\langle \hat{\Sigma},Q \right\rangle$, we can conclude that $Q^*=0$.
			
			Let's consider the case when $s^* > 0$. According to Lemma \ref{schur}, \eqref{qq} is equivalent to $Q \succeq \frac{1}{s} \hat{q} \hat{q}^{\top}$. Since $\hat{\Sigma}\succ0$, and the objective is to minimize $\left\langle \hat{\Sigma},Q \right\rangle$, we can deduce that $Q^*=\frac{1}{s} \hat{q} \hat{q}^{\top}$, $s^*=\arg\min\limits_{s>0}\frac{1}{s} \hat{q}^{\top}\hat{\Sigma}\hat{q}+\kappa_{1}s=\frac{\Vert\hat{\Sigma}^{\frac{1}{2}}\hat{q}\Vert}{\sqrt{\kappa_{1}}}$, and $\left\langle \hat{\Sigma},Q^* \right\rangle+\kappa_{1}s^*=\sqrt{\kappa_1}\left\|\hat{\Sigma}^{1 / 2}\hat{q}\right\|$.
			
			For the above two cases, after substituting $q=-2\left(\hat{q}+\Lambda \hat{\mu}\right)$, the Lagrange dual \eqref{qqobj}--\eqref{qq} of \eqref{max1} simplifies to
			\begin{eqnarray*}
				\min\limits_{r, \Lambda, Q, q, s}&& r+\left\langle \kappa_{2}\hat{\Sigma},\Lambda\right\rangle+(\Lambda^{\top}\hat{\mu}+q)^{\top}\hat{\mu}+ \tau_{1}\|x\|^{2}+\tau_{2}\alpha+\sqrt{\kappa_1}\left\|\hat{\Sigma}^{1 / 2}\left(q+2 \Lambda \hat{\mu}\right)\right\|\\
				\text { s.t.}\quad&&\psi(\xi_{a} -\xi_{B}^{\top}x)-(\Lambda^{\top}\xi+q)^{\top}\xi+\frac{\tau_{2}}{1-\beta}\left[-x^{\top}\xi_{B}-\alpha\right]_{+}-r\leq 0,~\forall \xi \in \mathcal{S}, \nonumber\\
				&& r\in \mathbb{R},\quad q\in\mathbb{R}^{d+1}, \quad\Lambda\in\mathbb{S}_+^{d+1},
			\end{eqnarray*}
			which is \eqref{innerdual}.

			Following Shapiro's duality theory for moment problems (Proposition 3.4 of \cite{shapiro2001duality}, (2.3) of \cite{xu2018distributionally}), the Slater type condition of \eqref{max1} can be written as later
			\begin{eqnarray}\label{slater1}
				(1,0,0)\in \text{int}\left\lbrace\left( \left\langle P,1\right\rangle,\left\langle P,\left( \begin{aligned}
					&-\hat{\Sigma}&\hat{\mu}-\xi\\
					&(\hat{\mu}-\xi)^\top &-\kappa_{1} 
				\end{aligned}\right)\right\rangle,\left\langle P,(\xi-\hat{\mu})(\xi-\hat{\mu})^{\top}-\kappa_{2} \hat{\Sigma}\right\rangle \right) \right. \nonumber\\
				\left. \quad+\{0\}\times\mathbb{S}_+^{d+2}\times \mathbb{S}_+^{d+1}\vert P\in\mathcal{M}\right\rbrace,
			\end{eqnarray}
			where $\left\langle P,\Psi(\xi)\right\rangle=\int_{\Xi}\Psi(\xi)P(d\xi)$. By Example 2.3 of \cite{xu2018distributionally}, the moment constraints \eqref{ambiguityset} satisfy the Slater type condition \eqref{slater1}. Then the equivalence between problems \eqref{max1} and \eqref{innerdual} holds. {Since the support set $\Xi$ is compact, we know that the optimal value of \eqref{max1} is finite. According to Proposition 3.4 of \cite{shapiro2001duality}, if the common optimal value of the primal problem and the dual problem is finite, then the set of optimal solutions of the dual problem is nonempty. Consequently, we can deduce that the set of optimal solutions of the dual problem \eqref{innerdual} is nonempty and bounded.}
		\end{proof}
		
		
		\section{Equivalent tractable reformulation of DRCVaR model}\label{equivalent}

		In this section, we establish an equivalent tractable reformulation for DRCVaR using Lagrangian dual. With the definition of the CVaR \eqref{cvar}, the problem \eqref{orip} is a ``min-max-min" optimization problem. Below, we first show that \eqref{orip} can be simplified as a nonsmooth SDP problem. Our derivation partially utilizes a similar argument for the Lagrange dual as in Delage et al. \cite{delage2010distributionally}. In addition to Delage's dual formulation, we exploit the nonsmooth function's structure and obtain a nonsmooth SDP. {For sake of simplicity, we denote 
			$$\nu = (x,\alpha,q,\Lambda),\quad\mbox{and}\quad \mathcal{V}=\Delta_{d}\times\mathbb{R}\times\mathbb{R}^{d+1}\times\mathds{S}_+^{d+1}. $$
			Below we will reformulate \eqref{orip} to the following  nonsmooth SDP 
			\begin{eqnarray}\label{pelconti}
				\min\limits_{\nu\in\mathcal{V}}~ \{\varphi(\nu):=\max\limits_{\xi\in\Xi}~ h_{\xi}(\nu)\},
			\end{eqnarray}
			where 
			\begin{eqnarray}\label{h}
				h_{\xi}(\nu):=h_{1}(\nu)+\sqrt{\kappa_{1}}\Vert\hat{\Sigma}^{\frac{1}{2}}(q+2\Lambda\hat{\mu})\Vert+h_{2,\xi}(x,q,\Lambda),
			\end{eqnarray}
			and $h_1$, $h_{2,\xi}$ are defined in \eqref{h1}.
			
			\begin{definition}[Section 2.2 of \cite{xu2018distributionally}]\label{weaklycompact}\\
				(i) If $\Xi$ is a compact set, then the set of all probability measures on $(\Xi, \mathds{B})$ is weakly compact with respect to (w.r.t.) topology of weak convergence.\\
				(ii) For a set of probability measures $\mathcal{\overline{A}} \subset \mathcal{M}$, $\mathcal{\overline{A}}$ is said to be weakly compact w.r.t. topology of weak convergence if every sequence $\left\{P_{N}\right\} \subset \mathcal{\overline{A}}$ contains a subsequence $\left\{P_{N^{\prime}}\right\}$, and a probability measure $P \in \mathcal{\overline{A}}$ such that $P_{N^{\prime}} \rightarrow P$.\\
				(iii) Under the topology of weak convergence, $\mathcal{\overline{A}}$ is said to be closed if for any sequence $\left\{P_{N}\right\}\subset \mathcal{\overline{A}}$ with $P_{N}\rightarrow P$ weakly, we have $P\in \mathcal{\overline{A}}$.
			\end{definition}
			Since $\Xi$ is a compact set in \eqref{ambiguityset}, by Definition \ref{weaklycompact} and the boundedness of probability measures on a compact set, the set of all probability measures on $(\Xi, \mathds{B})$ is compact under the topology of weak convergence.
			
			\begin{lemma}[Theorem 4.2 of \cite{sion1958general}]\label{minmax}
				Let $\hat{A}$ be compact, $\hat{B}$ be any space, $\rho(P,\alpha)$ be a function that is concave in $\hat{A}$ w.r.t. $P$, and convex in $\hat{B}$ w.r.t. $\alpha$. If $\rho(P,\alpha)$ is upper semi-continuous in $P$ for each $\alpha$, then $$\max\limits_{P\in \hat{A}} \min\limits_{\alpha\in \hat{B}} \rho(P,\alpha)=\min\limits_{\alpha\in \hat{B}} \max\limits_{P\in \hat{A}} \rho(P,\alpha).$$
			\end{lemma}

			\begin{theorem}[Nonsmooth SDP reformulation of DRCVaR]
				Consider the DRCVaR problem \eqref{orip}, and its associated nonsmooth SDP problem \eqref{pelconti}. Then, we have
				$$
				\min \eqref{orip}=\min \eqref{pelconti},
				$$
				in the sense that $x$ is an optimal solution of \eqref{orip} if and only if there exist $\alpha\in \mathbb{R}$, $q\in\mathbb{R}^{d+1}$, and $\Lambda\in\mathbb{S}_+^{d+1}$ such that $(x,~\alpha,~q,~\Lambda)$ is an optimal solution of \eqref{pelconti}.
			\end{theorem}
			\begin{proof}
				{Based on \eqref{K}, the problem \eqref{orip} is equivalent to the following problem:
					\begin{eqnarray}\label{mmm}
						\min\limits_{x\in \Delta_{d}}\max\limits_{P\in\mathscr{P}}\min\limits_{\alpha\in \mathbb{R}}~\int_{\Xi}\hat{K}(x,\alpha,\xi)P(d\xi).
					\end{eqnarray}
					Given that the support set $\Xi$ of the ambiguity set \eqref{ambiguityset} is compact, as per Definition \ref{weaklycompact}, we can establish that Lemma \ref{minmax} is applicable to the $\max\limits_{P\in\mathscr{P}}\min\limits_{\alpha\in \mathbb{R}}$ expression in \eqref{mmm}. This is due to the convexity of the function $\hat{K}(x,\alpha,P)$ with respect to $\alpha$ (stemming from the convexity of linear and plus functions) and the concavity (indeed, linearity) with respect to $P$, along with the compactness of $\mathscr{P}$. Then, interchanging the $\max\limits_{P\in\mathscr{P}}$ and $\min\limits_{\alpha\in \mathbb{R}}$ operators in \eqref{mmm} leads to an equivalent transformation
					\begin{eqnarray}\label{mimim}
						\min\limits_{x\in \Delta_{d}}\min\limits_{\alpha\in \mathbb{R}}\max\limits_{P\in\mathscr{P}}~\int_{\Xi}\hat{K}(x,\alpha,\xi)P(d\xi).
					\end{eqnarray}
					We start by considering the question of solving the inner maximization problem of \eqref{mimim} that uses the ambiguity set \eqref{ambiguityset}.}

				{Given any fixed $x$ and $\alpha$, we get the equivalence between \eqref{max1} and its dual problem \eqref{innerdual} by Proposition \ref{propinn}. We now return to discuss the semi-infinite program \eqref{innerdual}.} The main difficulty in solving a semi-infinite program comes from there being infinitely many constraints since infinitely many $\xi$'s values in the sample space $\Xi$. 
				We rewrite $$h_{2,\xi}(x,q,\Lambda)\leq r ,~~\forall\xi\in\Xi$$ into
				\begin{eqnarray*}
					\max\limits_{\xi\in\Xi}\left\{h_{2,\xi}(x,q,\Lambda)\right\}\leq r.
				\end{eqnarray*}
				Then the Lagrange dual problem \eqref{innerdual} is equivalent to
				\begin{eqnarray*}
					\begin{aligned}
						\min\limits_{q,~\Lambda}\quad&\max\limits_{\xi\in\Xi}\quad h_{\xi}(\nu)\\
						\rm{s.t.}\quad& q\in\mathbb{R}^{d+1},\quad\Lambda\in\mathbb{S}_+^{d+1}.
					\end{aligned}
				\end{eqnarray*}
				We use the fact that min-min operators can be performed jointly. This leads to an equivalent formulation \eqref{pelconti} for DRCVaR.
			\end{proof}

			{So far, we have transformed the DRCVaR index tracking model \eqref{orip} equivalently reformulated it as the SDP in \eqref{pelconti} shown in Theorem 1. 
				Since $X$ is bounded and $\cal P$ is bounded, and the objective function in  \eqref{orip} is continuous, we know that the solution set of the original problem  \eqref{orip} is nonempty and bounded.  Using Theorem 1, the existence of the optimal solutions to the SDP \eqref{pelconti} is guaranteed.}

				\section{The discretization scheme}\label{disscheme}

				For the purpose of computation, we give a discretization model for the equivalent tractable reformulation. We show the existence of solutions for the discretization model under mild conditions. Convergence results of the optimal values and solutions of the discretization scheme to those of the original equivalent reformulation are also given under mild assumptions. 

				In what follows, we consider the discrete approximation of \eqref{pelconti}. Our first step is to develop a discretized approximation of the continuous support set $\Xi$. Let $\Xi_{[N]}=\{\xi^1,\cdots, \xi^{N}\}$ be independent and identically distributed (iid) samples of $\xi$ drawn by Monte Carlo sampling from the set $\Xi$. We consider the following discretization scheme of \eqref{pelconti}
				\begin{eqnarray}\label{peldis}
					\begin{aligned}
						\min\limits_{\nu\in\mathcal{V}}~  \{\varphi^N(\nu):=\max\limits_{\xi\in\Xi_{[N]}}~ h_{\xi}(\nu)\},
					\end{aligned}
				\end{eqnarray}
				where $h_{\xi}$ is defined in \eqref{h}. 
					Then the corresponding approximation to $\phi_\beta(x)$ in \eqref{cvar} is
					\begin{eqnarray}\label{cvarN}
						\phi_\beta^N(x)= \min\limits_{\alpha\in\mathbb{R}} \left\{\alpha+\frac{1}{(1-\beta)N}\sum\limits_{i=1}^N\left[-x^{\top}\xi_{B}^i-\alpha\right]_{+}\right\}.
					\end{eqnarray}
					It is well known that the maximum w.r.t. $\alpha$ in the above formulation is achieved at a finite $\alpha$. In other words, we may restrict the maximum w.r.t. $\alpha$ to be taken within a closed interval $[-c, c]$ for some sufficiently large positive number $c$, see \cite{xu2018distributionally,rockafellar2000optimization}. Let $\mathbb{A}$ be the compact set consisting of the values of $\alpha$ for which the minimum in $\phi_\beta^N(x)$ is attained. 

					{\begin{proposition}\label{prop2.2}
							{The nonsmooth model in \eqref{peldis} is convex} and the solution set of \eqref{peldis} is nonempty.
					\end{proposition}}
					\begin{proof}
						{According to \eqref{h1} and \eqref{h}, we can see that $h_{\xi}(\nu)=h_{\xi}(x,\alpha,q,\Lambda)$, in which all the terms involving $x$, $\alpha$, $q$, and $\Lambda$ are block separable and convex with respect to the corresponding block variable, except that the term $$\sqrt{\kappa_{1}}\Vert\hat{\Sigma}^{\frac{1}{2}}(q+2\Lambda\hat{\mu})\Vert$$ involves both $q$ and $\Lambda$. 
							Fortunately, for any $\lambda \in (0,1)$, $q^{1}, q^{(2)} \in \mathbb{R}^{d+1}$,  and $\Lambda^{(1)}, \Lambda^{(2)} \in \mathds{S}_+^{d+1}$, we have by direct computation and the convexity of $\|\cdot\|$ that
							\begin{eqnarray*}
								& \quad \sqrt{\kappa_{1}}\left\Vert\hat{\Sigma}^{\frac{1}{2}}\left( (\lambda q^{(1)} + (1-\lambda) q^{(2)})+2(\lambda_1 \Lambda^{(1)} + (1-\lambda) \Lambda^{(2)})\hat{\mu}\right)\right\Vert \\
								&\le   \lambda\sqrt{\kappa_{1}}\left\Vert\hat{\Sigma}^{\frac{1}{2}}(q^{1}+2\Lambda^{(1)}\hat{\mu})\Vert + (1-\lambda)\sqrt{\kappa_{1}}\Vert\hat{\Sigma}^{\frac{1}{2}}(q^{(2)}+2\Lambda^{(2)}\hat{\mu})\right\Vert.
							\end{eqnarray*}
							Therefore, using the definition of convex function, we can show without difficulty that $h_{\xi}(\nu)$ with respect to $\nu$ is a convex function. By proposition 1.38 of \cite{mordukhovich2014easy}, the maximum function of a finite number of convex functions is also convex. In other words, $\varphi^N(\nu)$ is convex in $\nu$. Moreover, the feasible set $\mathcal{V}$ is a convex set. Hence, the nonsmooth model \eqref{peldis} is convex.}
						
						{From \eqref{cvarN}, we get if the sampling generates a collection of vectors $\xi_{B}^1,\ldots,\xi_{B}^N$, then the maximum w.r.t. $\alpha$ in \eqref{cvarN} is achieved at a finite $\alpha$. In other words, we may restrict the maximum w.r.t. $\alpha$ to be taken within a closed interval $[-c, c]$ for some sufficiently large positive number $c$, see \cite{xu2018distributionally,rockafellar2000optimization}. Let $\mathbb{A}$ be the compact set consisting of the values of $\alpha$ for which the minimum in $\phi_\beta^N(x)$ is attained.} By the compactness of $\Delta_{d}$ and $\Xi_{[N]}$, there exists $x^*\in\Delta_{d}$ and $\xi^*\in\Xi_{[N]}$ as the optimal solution of \eqref{peldis}. Then, the solution set of \eqref{peldis} is nonempty.
					\end{proof}
					
					{We respectively denote the optimal values of \eqref{peldis} and \eqref{pelconti} as $\hat{\vartheta}_{N}$ and $\hat{\vartheta}$}. 
					In Theorem \ref{the4.3}, we state convergence of problem \eqref{peldis} to problem \eqref{pelconti} in terms of the optimal value. 

					\begin{theorem}\label{the4.3}
						Suppose that the optimal value of \eqref{max1} is finite, and $\xi^1,\cdots, \xi^{N}$ be iid samples of $\xi$ and attaching to it with a continuous probability distribution $P$ over $\Xi$ such that
						\begin{eqnarray}\label{condixi}
							P\left(\left\|\xi-\xi^0\right\|<\delta\right) \geq C_2 \delta^{\gamma_2}
						\end{eqnarray}
						for any fixed point $\xi^0 \in \Xi$ and $\delta \in\left(0, \delta_0\right)$, where $C_2$, $\gamma_2$ and $\delta_0$ are some positive constants. When $N$ is sufficient large, for any positive number $\varepsilon$, there exist positive constants $\hat{C}(\varepsilon)$ and $\hat{\beta}(\varepsilon)$ such that
						\begin{eqnarray}\label{valcon}
							\mathrm{Prob}\left( \vert \hat{\vartheta}_{N}-\hat{\vartheta}\vert\geq\varepsilon\right)\leq \hat{C}(\varepsilon)e^{-\hat{\beta}(\varepsilon)N}.
						\end{eqnarray}
				\end{theorem}
				\begin{proof}
					(i) Denote by $$M_{\nu}(t):=\mathbb{E}\left[e^{t(h(\nu, \xi)-\mathbb{E}[h(\nu, \xi)])}\right]$$
					the moment generating function of the random variable $h(\nu, \xi)-\mathbb{E}[h(\nu, \xi)]$. Since $\Xi$ is a compact set and Section 3.1 of \cite{xu2018distributionally}, for each $\nu\in\mathcal{V}$, $\sup\limits_{\xi \in \Xi} h(\nu, \xi)<\infty$ and the moment generating function $M_{\nu}(t)$ is finite valued for all $t$ in a neighborhood of zero. 
					
					(ii) Through the continuity of the function $h$, we can establish the existence of a nonnegative measurable function $\kappa: \Xi \rightarrow \mathbb{R}_{+}$ and a constant $\gamma>0$ such that for all  $\xi \ in  \Xi$:
					\begin{eqnarray*}
						\left\vert h\left(\nu^{\prime}, \xi\right)-h\left(\nu^{\prime \prime}, \xi\right) \right\vert \leq \kappa(\xi)\left\|\nu^{\prime}-\nu^{\prime \prime}\right\|^\gamma, \quad\forall \nu^{\prime}, \nu^{\prime \prime} \in \mathcal{V}.
					\end{eqnarray*}
					
					(iii) Considering the boundedness of the support set $\Xi$ and Section 5 of \cite{shapiro2008stochastic}, we can conclude that the moment generating function $M_\kappa(t)$ of $\kappa(\xi)$ is finite for all $t$ in a neighborhood of zero.
					
					Combining the conditions mentioned above (i)--(iii), along with \eqref{condixi}, and the continuity of $h(\nu, \cdot)$ over $\Xi$, we can conclude that the relationship between the optimal values of \eqref{peldis} and \eqref{pelconti}, i.e., \eqref{valcon} holds true, as indicated by Lemma 3.1 (i) in \cite{xu2018distributionally}.
				\end{proof}
				
				\begin{remark}
					By Proposition 1 of \cite{anderson2020varying}, the condition \eqref{condixi} is very weak: it can be guaranteed whenever the density function of $\xi$ is lower bounded by a positive real-valued function which is analytic. It is easily seen to hold when the density function is bounded away from zero around $\xi^0$, and so it is effectively a condition on the way that the density approaches zero at the boundary of its support.
				\end{remark}
				
				\begin{remark}
					By Example 3 and Theorem 4 in \cite{xu2018distributionally}, we know that \eqref{orip} can be reformulated as an SDP, which can be solved by using an interior point (IP) method as discussed in \cite{delage2010distributionally}. But it requires the support set $\Xi$ of the random variable $\xi$ to be a compact ellipsoidal set. In this paper, we do not impose any restrictions on the form of the support set $\Xi$.
					
					It is worth mentioning that in \cite{xu2018distributionally}, the functions within the maximization operator are restricted to be differentiable with respect to $\nu$. However, in our model \eqref{peldis}, $h_{\xi}(\nu)$ incorporates non-differentiable terms arising from the $\ell_2$-norm and the plus function. Consequently, Algorithm 3.1 proposed in \cite{xu2018distributionally} is not directly applicable to solving our model \eqref{peldis}. Below we construct a smoothing function and employ the SPG method to solve our model \eqref{peldis}. 
				\end{remark}

				\section{SPG method}\label{sec4}

			\subsection{Smoothing function}\label{Stech}
			
			\begin{definition}(Definition 1 of  \cite{chen2012smoothing}) \label{smoothingdefinition}
				Let $g: \mathcal{V} \rightarrow \mathbb{R}$ be a continuous function. We call $\tilde{g}:\mathcal{V}\times \mathbb{R}_+ \rightarrow \mathbb{R}$ a smoothing function of $g$, if $\tilde{g}_{\mu}(\cdot)$ is continuously differentiable in $\mathcal{V}$ for any fixed $\mu > 0$, and for any $\nu \in \mathcal{V}$,
				\begin{eqnarray*}
					\lim\limits_{z\rightarrow \nu,~\mu \downarrow 0}\tilde{g}_{\mu}(z) = g(\nu)
				\end{eqnarray*}
				and $\left\{\lim\limits_{z\rightarrow \nu,~ \mu\downarrow 0}\nabla_{\nu}\tilde{g}_{\mu}(z)\right\}$ is nonempty and bounded.
			\end{definition}

			{\begin{definition}(Definition 2.3.4 of \cite{clarke1990optimization}) \label{regulardefinition}
					A function $g$ is said to be regular at $w$ provided
					\\
					(i) For all $v$, the usual one-sided directional derivative $f'(w;v)$ exists,
					\\
					(ii) For all $v$, $g'(w;v)=g^\circ(w;v)$, where
					\begin{eqnarray*}
						g'(w;v):=\lim\limits_{
							\lambda\downarrow 0}\frac{g(y+\lambda v)-g(y)}{\lambda}, ~g^\circ(w;v):=\limsup\limits_{
							y\rightarrow \nu,~
							\lambda\downarrow 0}\frac{g(y+\lambda v)-g(y)}{\lambda}.
					\end{eqnarray*}
			\end{definition}}
			
			\begin{lemma}(Proposition 1 of \cite{chen2012smoothing})\label{lh2} 
				Let $\tilde{\Gamma}_\mu$ and $\tilde{\pi}_\mu$ be smoothing functions of locally Lipschitz functions $\Gamma: \mathbb{R}^{N} \rightarrow \mathbb{R}$ and $\pi: \mathcal{V}\times\Xi_{[N]} \rightarrow \mathbb{R}^{N}$, respectively. The vector function
				\begin{eqnarray*}
					W(\nu,\Xi_{[N]})=\left(\pi(\nu,\xi^1), \ldots, \pi(\nu,\xi^N)\right)^T.
				\end{eqnarray*}
				If $\Gamma$ is regular at $W(\nu,\Xi_{[N]})$, $\pi$ is regular at $(\nu,\xi)$ and $\left.\nabla_z \tilde{\Gamma}_{\mu}(z)\right|_{z=W(\nu,\Xi_{[N]})} \geq 0$. Then $\tilde{\Gamma}_{\mu}(\tilde{W}_{\mu})$ is a smoothing function of $\Gamma(W)$ for any $\nu \in \mathcal{V}$.
			\end{lemma}
			
			\begin{lemma}(Proposition 2.5.6 of \cite{clarke1990optimization})\label{regular}
				Let $\hat{g}$ be Lipschitz near $\nu$. If $\hat{g}$ is convex, then $\hat{g}$ is regular at $\nu$.
			\end{lemma}

			Now we give the smoothing composite function using Lemmas \ref{lh2} and \ref{regular} for the following nonsmooth function $\Phi(H)$. For a vector function 
			\begin{eqnarray*}
				H(\nu,\Xi_{[N]})=\left(h(\nu,\xi^1), \ldots, h(\nu,\xi^N)\right)^T
			\end{eqnarray*}
			with components $h(\nu,\xi^i): \mathcal{V}\times\Xi \rightarrow \mathbb{R}$, {it is clear that} 
			\begin{eqnarray}\label{psi}
				{\varphi^N(\nu)}=\Phi(H(\nu,\Xi_{[N]}))=\max\{H(\nu,\Xi_{[N]})\}=\max\{h(\nu,\xi^1), \ldots, h(\nu,\xi^N)\}.
			\end{eqnarray}
			
			\begin{proposition}\label{prop3.2}
				For any fixed $\xi$ and the smoothing parameter $\mu>0$, the smoothing functions of $h(\nu,\xi)$ in \eqref{h} and $\varphi^N(\nu)$ in \eqref{peldis} respectively defined by
				\begin{eqnarray*}
					\tilde{h}_{\mu}(\nu,\xi)&=&h_{1}(\nu)+\sqrt{\kappa_{1}(q+2\Lambda\hat{\mu})\hat{\Sigma}(q+2\Lambda\hat{\mu})+\mu}+\tilde{\psi}_{\mu}(\xi_{a} -\xi_{B}^{\top}x)\\
					&~&\quad-(\Lambda^{\top}\xi+q)^{\top}\xi+\frac{\tau_{2} \mu}{1-\beta}\ln\left( 1+e^{\frac{-x^{\top}\xi_{B}-\alpha}{\mu}}\right) ,
				\end{eqnarray*}
				and
				\begin{eqnarray}\label{smoothingfun}
					{\tilde{\varphi}_{\mu}^N (\nu)}=\tilde{\Phi}_{\mu}(\tilde{H}_{\mu}(\nu,\xi_{[N]}))=\mu\ln \sum\limits_{i=1}^{N}e^{\frac{\tilde{h}_{\mu}(\nu,\xi^i)}{\mu}}.
				\end{eqnarray}
				In this context, if $\psi(c) = \|c\|_1$, then $\tilde{\psi}_{\mu}(c)=\sqrt{c^2+\mu}$. On the other hand, if $\psi(c) = \|c\|^2$, then $\tilde{\psi}_{\mu}(c)=\psi(c)$ for all $\mu\ge 0$.
			\end{proposition}
			\begin{proof}
				According to the smoothing technique in Example 2 of \cite{chen2012smoothing}, we have $\mu\ln(1+e^{\frac{t}{\mu}})$ is a smoothing function of $[t]_{+}$. Employing it into \eqref{h}, we have a smoothing function $\tilde{h}_{\mu}(\nu,\xi)$ of $h(\nu,\xi)$. Since the convexity of $h(\nu, \xi)$, and Lemma \ref{regular}, $h(\nu, \xi)$ is regular at $(\nu, \xi)$. Meanwhile, by Example 2 of \cite{chen2012smoothing}, we get the smoothing function of $\Phi(\hat{\chi})=\max\{\hat{\chi}_1,\ldots,\hat{\chi}_N\}$, which is $\tilde{\Phi}_\mu(\hat{\chi})=\mu\ln \sum\limits_{i=1}^{N}e^{\frac{\hat{\chi}_i}{\mu}}$ with 
				\begin{eqnarray*}
					\left.\nabla_{\hat{\chi}} \tilde{\Phi}_\mu(\hat{\chi})\right|_{\hat{\chi}=H(\nu,\Xi_{[N]})} =
					\left(\frac{e^{\frac{\tilde{h}_{\mu}(\nu,\xi^1)}{\mu}}}{\sum\limits_{i=1}^{N}e^{\frac{\tilde{h}_{\mu}(\nu,\xi^i)}{\mu}}},~\ldots,~\frac{e^{\frac{\tilde{h}_{\mu}(\nu,\xi^N)}{\mu}}}{\sum\limits_{i=1}^{N}e^{\frac{\tilde{h}_{\mu}(\nu,\xi^i)}{\mu}}}\right)^\top> 0.
				\end{eqnarray*}
				By the convexity of $h(\nu,\xi)$ and $\Phi(\hat{\chi})$, $\left.\nabla_{\hat{\chi}} \tilde{\Phi}_\mu(\hat{\chi})\right|_{\hat{\chi}=H(\nu,\Xi_{[N]})}> 0$, Lemmas \ref{lh2} and \ref{regular}, $\tilde{\varphi}_{\mu}^N(\nu)$ in \eqref{smoothingfun} is a smoothing function of $\varphi^N(\nu)$.
			\end{proof}

			\subsection{Smoothing projected gradient (SPG) method}

				We employ the SPG method \cite{zhang2009smoothing,zhang2019robust} to solve \eqref{peldis}, i.e. the discretization scheme of the equivalent reformulation.
				
				\vspace{3mm}\hspace{-5mm}\textbf{Algorithm 1}~\quad SPG method\label{algorithm1}\vspace{1mm}
				
				Initialize $\alpha_0, \sigma, \rho, \omega \in(0,1), \mu_0, \eta>0, \epsilon>0, \nu^0 \in \mathcal{V}$, and positive integer $n_0>0$. For $k \geq 0$ :\\
				\vspace{1mm}
				\hspace{5mm}Step 1. If
				\begin{eqnarray*}
					\left\|P_{\mathcal{V}}\left[\nu^k-\nabla_\nu {\tilde{\varphi}_{\mu_k}^N}\left(\nu^k\right)\right]-\nu^k\right\|<\epsilon,
				\end{eqnarray*}
				let $\nu^{k+1}=\nu^k$, and go to Step 3. Otherwise, go to Step 2.\\
				\vspace{1mm}
				\hspace{5mm}Step 2. Let $y^{0, k}=\nu^k$. For $j \geq 0: y^{j, k}(\alpha)=P_{\mathcal{V}}\left[y^{j, k}-\alpha \nabla_\nu {\tilde{\varphi}_{\mu_k}^N}\left(y^{j, k}\right)\right]$, and $y^{j+1, k}=$ $y^{j, k}\left(\alpha_{j, k}\right)$, where the stepsize $\alpha_{j, k}=\alpha_0 \rho^{\gamma_{j, k}}$ is chosen by the Armijo line search. That is, $\gamma_{j, k}$ is the smallest integer so that
				\begin{eqnarray*}
					{\tilde{\varphi}_{\mu_k}^N}\left(y^{j+1, k}\right) \leq {\tilde{\varphi}_{\mu_k}^N}\left(y^{j, k}\right)+\sigma\left\langle\nabla_\nu {\tilde{\varphi}_{\mu_k}^N}\left(y^{j, k}\right), y^{j+1, k}-y^{j, k}\right\rangle .
				\end{eqnarray*}
				If $j \geq n_0$ and $\frac{\left\|y^{j+1, k}-y^{j, k}\right\|}{\alpha_{j, k}}<\eta \mu_k$, set $\nu^{k+1}=y^{j+1, k}$, and go to Step 3 .\\
				\vspace{1mm}
				\hspace{5mm}Step 3. Choose $\mu_{k+1} \leq \omega \mu_k$.\\
				\vspace{3mm}
				
				Here, $j \geq n_0$ was proposed by \cite{zhang2009smoothing} to judge whether to go to Step 3 and decrease the smoothing parameter $\mu_k$. Thus for each $\mu_k$, at least $n_0$ inner iterations are performed. This modification helps to find a global or better local minimizer by the smoothing method from a computational point of view.

				{
					For any fixed $\bar{\nu} \in \mathcal{V}$, denote the Clarke subdifferential as
					\begin{eqnarray*}
						\partial \varphi^N(\bar{\nu})=\text{conv}\left\{U \mid \nabla \varphi^N(z) \rightarrow U, \quad\varphi^N ~\mathrm{is ~differentiable~ at~}~z,\quad z \rightarrow \bar{\nu}\right\},
				\end{eqnarray*}}
				and the subdifferential associated with a smoothing function as
				\begin{eqnarray*}
					G_{{\tilde{\varphi}_{\mu}^N}}(\bar{\nu})=\text{conv}\left\{U\mid \nabla_\nu {\tilde{\varphi}_{\mu}^N}\left(\nu^r, \mu_r\right) \longrightarrow U,\quad \text{for}\quad\nu^{r}\rightarrow\nu,\quad\mu_r \downarrow 0\right\} .
				\end{eqnarray*}
				From Definition \ref{smoothingdefinition}, it is clear that $G_{{\tilde{\varphi}_{\mu}^N}}(\bar{\nu})$ is a nonempty and bounded set. 
				\begin{definition}
					We say that $\nu^*$ is a stationary point of \eqref{peldis} associated with a smoothing function ${\tilde{\varphi}_{\mu}^N}$, if there exists $U \in G_{{\tilde{\varphi}_{\mu}^N}}\left(\nu^*\right)$ such that
					\begin{eqnarray*}
						\left\langle U, \nu^*-z\right\rangle \leq 0 \quad \text { for all } z \in \mathcal{V} .
					\end{eqnarray*}
				\end{definition}
				
				\begin{definition}
					We say that $\nu^*$ is a Clarke stationary point of \eqref{peldis}, if there exists $U \in \partial \varphi^N(v^*)$ such that
					\begin{eqnarray*}
						\left\langle U, \nu^*-z\right\rangle \leq 0 \quad \text { for all } z \in \mathcal{V} .
					\end{eqnarray*}
				\end{definition}
				
				\begin{theorem}
					Any accumulation point $\nu^*$ of $\left\{\nu^k\right\}$ generated by Algorithm 1 with the smoothing function $\tilde{\varphi}^N_{\mu}$ is a global minimizer  of (\ref{peldis}).
				\end{theorem}
				\begin{proof}
					According to the definition of the subdifferential associated with a smoothing function, $G_{{\tilde{\varphi}_{\mu}^N}}(\bar{\nu})$, we can conclude that any accumulation point $\nu^*$ of the sequence $\left\{\nu^k\right\}$ generated by Algorithm 1 with the smoothing function ${\tilde{\varphi}_{\mu}^N}$ belongs to the set $G_{{\tilde{\varphi}_{\mu}^N}}(\bar{\nu})$. By Proposition 1 in \cite{chen2012smoothing}, we have
					\begin{eqnarray*}
						\partial \varphi^N(\bar{\nu}) = G_{{\tilde{\varphi}_{\mu}^N}}(\bar{\nu}).
					\end{eqnarray*}
					Therefore, we can conclude that any accumulation point $\nu^*$ of the sequence $\left\{\nu^k\right\}$ generated by Algorithm 1 with the smoothing function ${\tilde{\varphi}_{\mu}^N}$ belongs to the set $\partial \varphi^N(\bar{\nu})$. In other words, any accumulation point $\nu^*$ of $\left\{\nu^k\right\}$ generated by Algorithm 1 with the smoothing function ${\tilde{\varphi}_{\mu}^N}$ is a Clarke stationary point of 
					(\ref{peldis}). Since (\ref{peldis}) is a convex model as shown in Proposition \ref{prop2.2}, we know that any Clarke stationary point of (\ref{peldis}) is a global minimizer of (\ref{peldis}).
				\end{proof}

				\section{Numerical results}\label{sec5}

				We consider the portfolio optimization problem \eqref{peldis} where the investor makes an optimal decision using a historical return rate of $d=120$ stocks 
				between January 2008 and July 2023 from National Association of Securities Deal Automated Quotations (NASDAQ) index\footnotemark[2]\footnotetext[2]{https://cn.investing.com}, which contains 3921 samples. All experiments are performed in Windows 10 on an Intel Core 10 CPU at 3.70 GHZ with 64 GB of RAM, using MATLAB R2021a.

				To evaluate the portfolios in an out-of-sampling setting, we employ a rolling window approach with a window size of 3500 daily observations for the above data, that is, $N=3500$ in \eqref{peldis}. The rolling window approach used for the daily data is the same as that in section 4 of \cite{kremer2020sparse}. The first $\tau$ return observations are used to compute the optimal weight vector $\hat{x}_t$. The resulting portfolio is assumed to be held for the following 21 days. 
				We utilize the data from these 21 days as a testing vector $B_{t+1}\in \mathbb{R}^{120}$ to compute the out-of-sample performance. Here, $(B_{t+1})_{i}$ denotes the ratio of the closing price on the last day of the 21 days to the opening price on the first day of the 21 days for the $i$-th asset, $a_{t+1}\in \mathbb{R}$ is the observed corresponding random market index return. During this period, the investment based on portfolio $\hat{x}_t$ gains wealth by $B_{t+1}^{\top}\hat{x}_t$. After these 21 days, based on the $B_{t+1}$, the portfolio $\hat{x}_t$ is transited to $\hat{x}_{t+}$, whose component is calculated by $(\hat{x}_{t+})_{i}=\frac{(B_{t+1})_{i}(\hat{x}_{t})_{i}}{B_{t+1}^{\top}\hat{x}_t}$, as pointed out in section 2 of \cite{yang2018reversion}. In the next step, we roll the data window forward, discarding the oldest and including the most recent 21 observations to our training set. Solving the model \eqref{peldis} with the new training set, we obtain a solution $\hat{x}_{t+1}$, and then we rebalance $\hat{x}_{t+}$ to the new weight vector $\hat{x}_{t+1}$, which determines our portfolio holdings and the out-of-sample return for the next 21 days. This process is repeated until the window can not be rolled anymore, i.e., $t=1,~\ldots,~\bar{t}=\left\lfloor\frac{N_{\rm{tol}}-\tau}{21}\right\rfloor$. In our setting, in view of $N_{\rm{tol}}=3921$ and $\tau=3500$, it follows that $\bar{t}=20$.
				
				The in-sample tracking error denoted by TEI, and the out-of-sample tracking error referred to as TEO will be taken into consideration to illustrate the performance of the given methods, which are defined by
				\begin{eqnarray*}
					{\rm TEI} = \frac{1}{\bar{t}} \sum_{t=1}^{\bar{t}}
					\left[ \frac{1}{\tau}\left\|(1,-\hat{x}_{t}^{\top})\left(\xi^{21(t-1)+2},\ldots,\xi^{21(t-1)+1+\tau} \right)\right\|^{2}\right] ,
				\end{eqnarray*}
				and
				\begin{eqnarray*}
					{\rm TEO} = \frac{1}{\bar{t}} \sum_{t=1}^{\bar{t}}
					\left[ \left\|(1,-\hat{x}_{t}^{\top})(a_{t+1};B_{t+1})\right\|^{2}\right] .
				\end{eqnarray*}
				A smaller TEI or TEO indicates better performance of the portfolio. We also evaluate the out-of-sample performance of the portfolio strategies using three measures: the portfolio variance ($\hat{\sigma}^2$), the Sharpe ratio (${\rm SR}$), and the turnover (${\rm TO}$), of each portfolio strategy:
				{\begin{eqnarray*}
						\hat{\sigma}^2=\frac{1}{\bar{t}-1} \sum_{t=1}^{\bar{t}}\left(B_{t+1}^{\top} \hat{x}_{t}-\hat{\mu}\right)^2, \quad \text { with } \ \hat{\mu}=\frac{1}{\bar{t}}\sum_{t=1}^{\bar{t}} B_{t+1}^{\top} \hat{x}_{t}, \\
						{\rm SR}=\frac{\hat{\mu}}{\hat{\sigma}}, \qquad\mbox{and}\qquad
						{\rm TO}=\frac{1}{\bar{t}} \sum_{t=1}^{\bar{t}}\left\|\hat{x}_{t+1}-\hat{x}_{t+}\right\|_1 .
				\end{eqnarray*}}
				A larger value of ${\rm SR}$ indicates better performance of a portfolio. While smaller values of TEI, TEO, $\hat{\sigma}^2$, and $\rm{TO}$ indicate better performance of the portfolio.
				
				For the parameters in the SPG method, we set
				\begin{eqnarray}\label{para}
					\alpha_0=1,~ \sigma=10^{-6},~ \rho=\frac{1}{2}, ~\mu_0=1, ~\eta=10^3, ~\omega=\frac{1}{2},~ \epsilon=10^{-4}, ~n_0=5 .
				\end{eqnarray}
				and we stop the iteration of the SPG method and set the solution to $\hat{\nu}_t=\nu^k$ if
				\begin{eqnarray}\label{stop}
					\left\|P_{\mathcal{V}}\left[\nu^k-\nabla_\nu {\tilde{\varphi}_{\mu_k}^N}\left(\nu^k\right)\right]-\nu^k\right\| \leq 10^{-4}, \quad \text { and } \quad \mu_k \leq 2 \times 10^{-6},
				\end{eqnarray}
				or the number of iterations exceeds 3000. 
				{We set $\kappa_{1}=0.1$ and $\kappa_{2}=1$ of \eqref{ambiguityset} that are the same as \cite{li2020fast} in our numerical experiments.}
				
				First, to see the effectiveness of the SPG method, we compare it with the subgradient (Subgrad) method in solving \eqref{peldis} with $\psi(a)=a^2$, using the first 3500 training samples. 
			The Armijo's step size rule is used in the Subgrad method. 
		We set $\tau_{1}=\tau_{2}=10^{-2}$, and plot in Figure \ref{fig1.1} the objective value in the training set (Obj) versus (v.s.) CPU time, of our SPG method and the Subgrad method. 

		\begin{figure}[htp]
			\center
			\includegraphics[height=5.5cm]{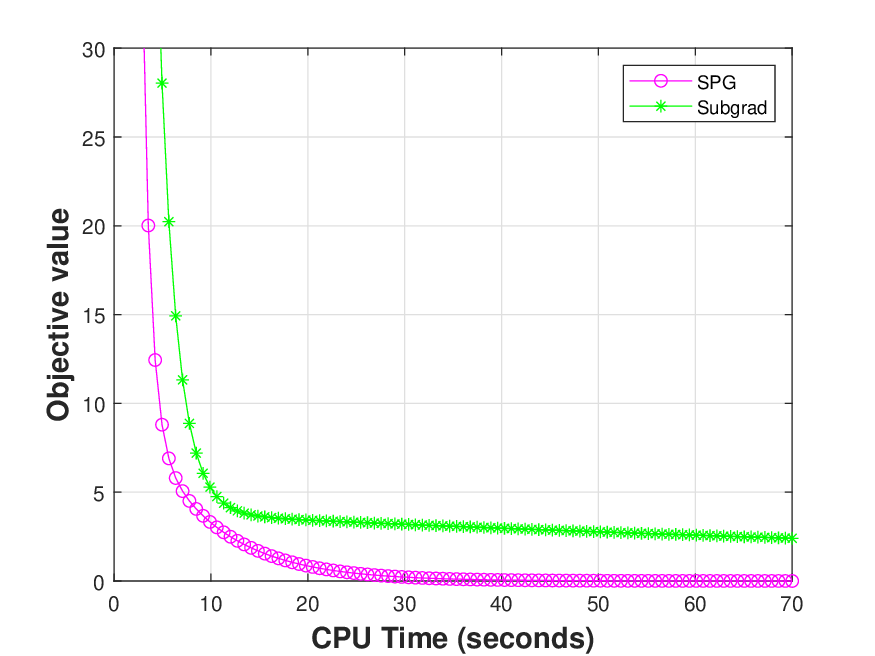}
			\caption{Obj v.s. CPU time, when $\tau_{1}=\tau_{2}=10^{-2}$}\label{fig1.1}
		\end{figure}
	\begin{table}[htp]
			\caption{In-sample and out-of-sample performances, and CPU time with  $\tau_{1}=\tau_{2}=10^{-2}$ }\label{tab1.1}
				\begin{footnotesize}
					\scalebox{1}{
						\begin{tabular}{lccccccc}
							\toprule
							\multirow{2}{*}{ALG.} &\multirow{2}{*}{$N$} &\multirow{2}{*}{$d$} &\multirow{2}{*}{$\mu$}	&\multicolumn{2}{c}{In-sample performance}&\multirow{2}{*}{TEO} &\multirow{2}{*}{CPU}  \\ \cline{5-6}
							& &	 & &Obj  &TEI &  &	\\ \midrule
							SPG &\multirow{2}{*}{3500}&\multirow{2}{*}{120}
							&$\downarrow$ &\textbf{0.4105}
							&\textbf{6.7795e-05}&\textbf{2.3298e-04 } &\textbf{69.9531}  \\
							Subgrad  & &
							&- &0.4249
							&1.4837e-04 &7.0829e-04&690.1406 \\
							\bottomrule
					\end{tabular}}
				\end{footnotesize}
		\end{table}

		Next, to see the merit of our DRCVaR index tracking model, we compare our model with the other index tracking models mentioned in this paper. The out-of-sample performances are evaluated. To be specific, we choose the optimal values of $\tau_{1},~\tau_{2}$  by varying them on the grid
		\begin{eqnarray}\label{grid}
			\{k \times 10^{-4},~k=0,~2,~4,~6,~8,~10\}
		\end{eqnarray}
		and the values with the lowest TEO are chosen for each model. We record in Table \ref{tab1.2} the lowest TEO, together with the corresponding parameters $\tau_{1},~\tau_{2}$, 
		in-sample and out-of-sample performances. {In Table \ref{tab1.2}, we give our DRCVaR model with $\ell_1$-norm $(\mathrm{DRCVaR}-\ell_1)$, and with $\ell_2$-norm $(\mathrm{DRCVaR}-\ell_2)$, respectively. The $(\mathrm{DRCVaR}-\ell_1)$ model outperforms the  $(\mathrm{DRCVaR}-\ell_1)$ model for all the evaluation criteria except TO. 
			
			Besides, we compare our DRCVaR model with all the other models in Section \ref{secmodel}, including $\mathrm{Lasso~sparse}$, $\mathrm{Mixed~0-1~LP}$, $\mathrm{TE-\ell_0}$, $\ell_2-\ell_{3/4}$, and $\ell_2-\ell_{1/2}$. Additionally, we also compare the stochastic models, that is, the $\mathrm{SCVaR-\hat{P}}$ model with $\ell_1$-norm $(\mathrm{SCVaR-\hat{P}}-\ell_1)$, and with $\ell_2$-norm $(\mathrm{SCVaR-\hat{P}}-\ell_2)$. Here, $\hat{P}$ refers to the empirical distribution obtained from historical data. 
			How to solve the above models is explained in Section \ref{sec3}. For parameters $\tau_1$ and $\tau_2$, ``-" indicates that the parameters are not required in the corresponding model.}
		\begin{table}[htp]
			\caption{The in-sample and out-of-sample performance, and parameters $\tau_{1},~\tau_{2}$ determined by the lowest TEO.}	\label{tab1.2}
			\begin{footnotesize}
				\scalebox{0.83}{	\begin{tabular}{lcccccccc}
					\toprule
					\multirow{2}{*}{Model} &\multirow{2}{*}{$\tau_{1}$}&\multirow{2}{*}{$\tau_{2}$}	&\multirow{2}{*}{TEI}  &\multicolumn{4}{c}{Out-of-sample performances}&\multirow{2}{*}{CPU}  \\ \cline{5-8}
					&		&  &  &TEO&{$\hat{\sigma}^2$}& ${\rm SR}$& ${\rm TO}$ &	\\ \midrule
					$\mathrm{DRCVaR}-\ell_1$	&2e-4 &8e-4  &6.7069e-05 &7.4505e-04&1.1206e-03 &14.8538  &\textbf{1.5711e-03} &89.1562 \\
					$\mathrm{DRCVaR}-\ell_2$	&2e-4 &6e-4  &\textbf{6.6353e-05} &\textbf{7.1391e-04}&\textbf{1.0115e-03} &\textbf{19.4181}  &1.5937e-03 &\textbf{86.4375} \\ 
					$\mathrm{Lasso~sparse}$ 	&4e-4 &-  &8.9501e-05 &2.5973e-03 &2.0911e-02&7.8261  &9.8089e-02 &313.2188 \\ 
					$\mathrm{Mixed~0-1~LP}$ 	&- &-  &1.2174e-04 &9.1707e-04 &1.0788e-03&13.5418 &2.2725e-03 &387.2490  \\ 
					$\mathrm{TE-\ell_0}$	&- &-  &1.3777e-04 &3.2221e-03 &5.1207e-03&13.1119 &2.7314e-03 &218.2344   \\ 
					$\ell_2-\ell_{3/4}$	&4e-4 &1e-3  &1.0742e-04 &8.6735e-03 &8.5709e-03&11.6775  &5.0517e-03 &95.6406   \\ 
					$\ell_2-\ell_{1/2}$	&2e-4 &8e-4  &1.1471e-04 &9.5441e-03 &9.6140e-03&11.0146 &7.0298e-03 &106.6388   \\ 
					$\mathrm{SCVaR-\hat{P}}-\ell_1$	&6e-4 &6e-4  &1.3212e-04 &2.5708e-03 &1.1925e-02&11.9054  &7.9472e-02 &117.1719   \\ 
					$\mathrm{SCVaR-\hat{P}}-\ell_2$	&6e-4 &1e-3  &1.1505e-04 &2.4398e-03 &1.0878e-02&12.5736  &8.1754e-02 &113.9688  \\ 
					\bottomrule
				\end{tabular}}
			\end{footnotesize}
		\end{table}
		Clearly, our $\mathrm{DRCVaR}$ models, when equipped with appropriate parameters and solved by the SPG method, exhibit strong performance both in-sample and out-of-sample. 
		Comparing the results of $\mathrm{DRCVaR}-\ell_2$, $\mathrm{DRCVaR}-\ell_1$ with $\mathrm{SCVaR-\hat{P}}-\ell_2$, and $\mathrm{SCVaR-\hat{P}}-\ell_1$, we can observe that DRO can enhance the out-of-sample performance of the models. Furthermore, the shorter CPU time of $\mathrm{DRCVaR}-\ell_2$, $\mathrm{DRCVaR}-\ell_1$ compared to $\mathrm{SCVaR-\hat{P}}-\ell_2$, and $\mathrm{SCVaR-\hat{P}}-\ell_1$ is attributed to the utilization of the SPG method.

							\section{Conclusion}\label{sec6}{
								In this paper, we propose a novel index tracking model that integrates the idea of the distributionally robust optimization and the conditional value-at-risk penalty. The ambiguity set is defined by the first two moments. We transform it into an equivalent nonsmooth convex minimization, and give an approximate discretization scheme for the continuous random vector. The objective function in this scenario encompasses the maximum of a finite number of nonsmooth functions. 
								We employ the SPG method to solve the discretized problem, Empirical results using the NASDAQ index dataset from January 2008 to July 2023 demonstrate the promising of our proposed model as well as the efficiency of the SPG method.}
\vskip 1mm

%


\begin{footnotesize}
\bibliographystyle{siamplain}
\bibliography{references}
\end{footnotesize}

\end{document}